\documentclass{article}

\usepackage{fancyhdr} 				% to manipulate headers and footers and pagestyles
\usepackage[fleqn]{amsmath} 	% to move the displaied equations to the left
\usepackage{amsthm,amsfonts,latexsym,amssymb,amscd} % AMS macros
\usepackage[mathscr]{eucal}   % to use the font eucal (\As ... \Zs below)  
\usepackage[all]{xy}					% for diagrams
\usepackage{stmaryrd} 				% for double parenthesis characters
\usepackage[draft=false]{hyperref} 	% to generate hyperlinks in pdf
\newcommand{\A}{\mathcal{A}}
\newcommand{\B}{\mathcal{B}}

\renewcommand{\H}{\mathcal{H}} 	%re!

\newcommand{\K}{\mathcal{K}}
\renewcommand{\L}{\mathcal{L}}	%re!

\renewcommand{\P}{\mathcal{P}}	%re!

\newcommand{\Fg}{\mathfrak{F}}

\newcommand{\Rg}{\mathfrak{R}}

\newcommand{\CC}{{\mathbb{C}}}

\newcommand{\NN}{{\mathbb{N}}}

\newcommand{\RR}{{\mathbb{R}}}
\newcommand{\TT}{{\mathbb{T}}}

\newcommand{\As}{{\mathscr{A}}}

\newcommand{\Fs}{{\mathscr{F}}}

\newcommand{\Rs}{{\mathscr{R}}}

\DeclareFontFamily{OT1}{rsfs}{}
\DeclareFontShape{OT1}{rsfs}{m}{n}{<10> <7> <5> gen * rsfs}{}
\DeclareMathAlphabet{\scr}{OT1}{rsfs}{m}{n} %\rs was substituted here with \scr 

\newcommand{\Af}{\scr{A}}
\newcommand{\Bf}{\scr{B}}

\newcommand{\Gf}{\scr{G}}
\newcommand{\Hf}{\scr{H}}

\newcommand{\Mf}{\scr{M}}

\newcommand{\Sf}{\scr{S}}

\DeclareMathOperator{\ke}{Ker}

\DeclareMathOperator{\spa}{span}
\DeclareMathOperator{\id}{Id}
\DeclareMathOperator{\Ob}{Ob}
\DeclareMathOperator{\Mor}{Mor}
\DeclareMathOperator{\car}{card}

\DeclareMathOperator{\Aut}{Aut}
\DeclareMathOperator{\ad}{ad}
\renewcommand{\emph}{\textbf} 									% emphasis with boldface
\newcommand{\cj}[1]{\overline{#1}}							% usage: \cj{x} for conjugate of x
\newcommand{\ip}[2]{\langle #1\mid #2\rangle}		% usage \ip{a}{b} Dirac Inner product
\renewcommand{\iff}{\Leftrightarrow}						% IFF
\newcommand{\imp}{\Rightarrow}									% implication
\newcommand{\cs}{$\text{C}^*$-algebra}					% C* algebra
\newcommand{\hlink}[1]{\href{#1}{\texttt{#1}}} % usage: \hlink{URL} to make hyperlink with text anchor given by their same URL in typewriter stile.

\newtheorem{theorem}{Theorem}[section]
\newtheorem{corollary}[theorem]{Corollary}
\newtheorem{lemma}[theorem]{Lemma}
\newtheorem{proposition}[theorem]{Proposition}

\newtheorem{definition}[theorem]{Definition}
\newtheorem{remark}[theorem]{Remark}

\title{A Category of Spectral Triples and Discrete Groups with Length Function.}

\author{Paolo Bertozzini$^*$\\
Department of Mathematics and Statistics\\ 
Faculty of Science and Technology\\ 
Thammasat University - Rangsit Campus\\
Bangkok 12121, Thailand\\
E-mail: paolo@mathstat.sci.tu.ac.th \\
{}\\
Roberto Conti$^*$ \\
Department of Mathematics \\
Faculty of Science \\
Chulalongkorn University \\
Bangkok 10330, Thailand \\
E-mail: 
conti@math.sc.chula.ac.th \\
{}\\
Wicharn Lewkeeratiyutkul$^*$ \\
Department of Mathematics \\
Faculty of Science \\
Chulalongkorn University \\
Bangkok 10330, Thailand \\
E-mail:	Wicharn.L@chula.ac.th
}
\date{February 01, 2005}

\numberwithin{equation}{section}  % to enumerate equations after sections        
\begin{document}

\maketitle

\begin{abstract} 
In the context of Connes' spectral triples, a suitable notion of morphism is introduced. Discrete groups with length function provide a natural example for our definitions. 
Connes' construction of spectral triples for group algebras is a covariant functor
from the category of discrete groups with  length functions to that of
spectral triples.
Several interesting lines for future study of the categorical properties
of spectral triples and their variants are suggested.

\medskip

\noindent 
\emph{MSC-2000:} 	46L87,			% Non-commutative Geometry 
									18F99, 			% Categories and Geometry
									20C07,				% Group Rings on Infinite Groups and Their Modules
									22D15.				% Group Algebras of Locally Compact Groups
					
\smallskip

\end{abstract}

\thanks{\noindent $^*$ Partially supported by the Thai Research Fund.}

\section{Introduction}

The notions of morphism, as a generalization of ``coordinate transformation'', and respectively of category, as a generalization of ``group of transformations'', are going to be central in all the attempts to reformulate the concepts of physical covariance in an algebraic context (see for instance J.~Baez~\cite{Ba}). 

In the abstract framework of A.~Connes' Non-commutative Geometry~\cite{C2,FGV}, where non-commutative manifolds are described by spectral triples, a definition of ``morphism of spectral triples'' is still missing in the literature. With the present short note, we intend to provide tentative definitions of ``morphism'' and of ``category of spectral triples'', and to investigate some of their properties.

Since, as typical feature of every non-commutative geometric setting, ``non-commutative spaces'' are described dually by the category of ``spectra'' (categories of representations) of their algebras of functions, defining a morphism of non-commutative spaces actually amounts to the specification of a functor between representations categories and, under this point of view, our work can also be seen as an example of ``categorification'' process in which sets are replaced by categories (see for example J.~Baez, J.~Dolan~\cite{BD} or L.~Ionescu~\cite{I}). 

In the second part of this paper, we proceed to the construction of a natural covariant functor, from the category of discrete groups equipped with a length function, to our category of spectral triples, that shows the validity of the proposed definition of morphisms. We expect this functor to be just one particular example in a class of functors from suitable categories of ``geometrical objects'' to the category of spectral triples. 

Actually this work is part of a much wider research project~\cite{BCL1} that, among several other objectives, has the purpose to study an appropriate notion of non-commutative (totally geodesic) submanifold and quotient manifold and the study of some suitable functorial relations between the categories of spectral triples and spin$^c$ Riemannian manifolds.  
This program will be carried out in detail in a forthcoming paper~\cite{BCL2}.
The situations investigated here are usuful to present all the relevant structures involved without dealing with the complications arising from ``spinorial calculus'' on Riemannian spin$^c$ manifolds.

\medskip

Treatments of non-commutative geometry in a suitable categorical framework, mostly appealing to Morita equivalence, have already appeared in a more or less explicit form. 
In~\cite{C3,C4,C5} A.~Connes shows how to transfer a given Dirac operator using Morita equivalence bimodules and compatible connections on them, thus leading to the concept of ``inner deformations'' of a spectral geometry that encompasses a formula for expressing the transformed Dirac operator in the form $\widetilde{D} = D + A + JAJ^{-1}$.
The categorical ``ideology'' becomes especially evident among the practitioners of ``non-commutative algebraic geometry'' (see for example M.~Kontsevich and A.~Rosenberg~\cite{KR1,KR2, R}) and morphisms between non-commutative manifolds, thought of as non-commutative spectra, have been proposed by Y.~Manin~\cite{Ma} in terms of the notion of ``Morita morphisms'', i.e.~functors among representations categories that are obtained by tensorization with bi-modules.

The notion of morphism of spectral triples described in the sequel is not as general as possible, and several further generalizations are undoubtedly at hand.

In a wider perspective~\cite{BCL3} a morphism of the spectral triples $(\As_j,\H_j,D_j)$, with $j=1,2$, might be formalized as a ``suitable'' functor $\Fs: {}_{\As_2}\Mf\to {}_{\As_1}\Mf$, between the categories ${}_{\As_j}\Mf$ of $\As_j$-modules, having ``appropriate intertwining'' properties with the Dirac operators $D_j$. 

The morphisms described in the sequel are only a very special case of that picture. However for the present purposes that level of generality would be unnecessary, and so we stick to the more restrictive definition provided by homomorphisms $\phi:\As_1\to \As_2$ of algebras with an intertwining operators $\Phi:\H_1\to\H_2$ between the Dirac operators.

We can thus establish our main result, stating that Connes' construction of spectral triples from group algebras is functorial in nature.

Whether these functors can be chosen to be full, if they are extendable to non-monomorphic cases and, in a broader context, which other functors into categories of spectral triples can be obtained this way seem to be interesting questions and we hope to return to these and related issues elsewhere.

\section{A Category of Spectral Triples.}

In this section we define a ``natural'' notion of morphism between spectral triples. 
Examples will be provided in the next section~\ref{sec: gr}.

In order to facilitate the reader and to establish our notations, we start recalling the  definitions and key properties related to spectral triples. 

\subsection{Preliminaries on Spectral Triples.}

\label{sec: sptr}

A.~Connes (see~\cite{C2, FGV}) has proposed a set of axioms for ``non-commutative manifolds''\footnote{At least in the case of compact, finite dimensional, Riemannian, orientable, spin$^c$ manifolds}, called a (compact) spectral triple or an (unbounded) K-cycle.
\begin{itemize}
\item
A (compact) \emph{spectral triple} $(\As, \H, D)$ is given by:
\begin{itemize}
\item
a unital pre-C$^*$-algebra\footnote{Sometimes $\As$ is required to be closed under holomorphic functional calculus.} $\As;$
\item
a representation $\pi: \As \to \B(\H)$ of $\As$ on the Hilbert space $\H;$
\item
a (non-necessarily bounded) self-adjoint operator $D$ on $\H$, called the Dirac operator, such that:
\begin{itemize}
\item[a)]
the resolvent $(D-\lambda)^{-1}$ is a compact operator, $\forall \lambda \in \CC\setminus\RR$,\footnote{As already noticed by Connes, 
this condition has to be weakened in the case of non-compact manifolds, cf. \cite{GLMV,GGISV,Re1,Re2}.}
\item[b)]
$[D,\pi(a)]_{-}\in \B(\H),$ 
for every $a \in \As,$ \\ 
where $[x,y]_{-}:= xy- yx$ denotes the commutator of $x,y \in \B(\H).$ 
\end{itemize}
\end{itemize}
\item
The spectral triple is called \emph{even} if there exists a grading operator, i.e.~a bounded self-adjoint operator $\Gamma \in \B(\H)$ such that:
\begin{gather*}
\Gamma^2=\text{Id}_\H; \quad  [\Gamma, \pi(a)]_{-}=0, \forall a \in \As; \quad
[\Gamma, D]_{+}=0, 
\end{gather*}
where $[x,y]_{+}:=xy+yx$ is the anticommutator of $x,y.$ 

A spectral triple that is not even is called \emph{odd}. 
\item
A spectral triple is \emph{regular} if
the function
$$\Xi_x: t\mapsto \exp(it|D|)x\exp(-it|D|)$$ 
is regular, 
i.e.~$\Xi_x\in \text{C}^\infty(\RR, \B(\H)),$\footnote{
This condition is equivalent to 
$\pi(a), [D,\pi(a)]_{-} \in \cap_{m=1}^\infty \text{Dom}\, \delta^m,$ for all $a\in \As$, where $\delta$ is the derivation given by $\delta(x):=[|D|, x]_{-}$.}
for every $x \in \Omega_D(\As)$, 
where~\footnote{We assume that for $n=0 \in \NN$ the term in the formula simply reduces to $\pi(a_0)$.}
\begin{equation*}
\Omega_D(\As):=\spa \{\pi(a_0)[D,\pi(a_1)]_- \cdots [D,\pi(a_n)]_- \  | 
\ \  n\in \NN, \ a_0, \dots, a_n \in \As\} \ .
\end{equation*}
\item
The spectral triple is \emph{$n$-dimensional} iff there exists an integer $n$ such that the Dixmier trace of $|D|^{-n}$ is finite nonzero.
\item
A spectral triple is \emph{$\theta$-summable} if $\exp(-tD^2)$ is a trace-class operator for every $t>0.$
\item
A spectral triple is \emph{real} if there exists an antiunitary operator $J: \H \to \H$ such that: 
\begin{gather*}
[\pi(a), J\pi(b^*)J^{-1}]_{-}=0, \quad \forall a,b \in \As; \\
[\, [D, \pi(a)]_{-}, J\pi(b^*)J^{-1}]_{-}=0, \quad \forall a,b \in \As, \quad {\emph{first order condition};} \\
J^2=\pm\text{Id}_\H;  \quad [J,D]_{\pm}=0; \\
\text{and, only in the even case,} \quad
[J,\Gamma]_{\pm}=0,
\end{gather*}
where the choice of $\pm$ in the last three formulas depends on the ``dimension'' $n$ of the spectral triple modulo $8$ in accordance to  the following table: 
\begin{center}\label{tb: J}
\begin{tabular}{|l|c|c|c|c|c|c|c|c|}
\hline
$n$									&$0$	&$1$	&$2$	&$3$	&$4$	&$5$	&$6$	&$7$	\\
	\hline
$J^2=\pm\text{Id}_\H$	&	$+$	&	$+$		&	$-$		&	$-$		&	$-$		&	$-$		&	$+$		& $+$		\\
	\hline
$[J,D]_{\pm}=0$				&	$-$	&	$+$		&	$-$		&	$-$		&	$-$		&	$+$		&	$-$		& $-$		\\
\hline
$[J,\Gamma]_{\pm}=0$	&	$-$	&			&	$+$		&			&	$-$		&			&	$+$		& 		\\
\hline
\end{tabular}
\end{center}
\item
A spectral triple is called \emph{commutative} if the algebra $\As$ is commutative.
\end{itemize}

\subsection{Morphisms of Spectral Triples.}

The objects of our category $\Sf$ will be spectral triples $(\As,\H,D).$
Given two spectral triples $(\As_j,\H_j,D_j),$ with $j=1,2,$ a \emph{morphism of spectral triples} is a pair 
\begin{gather*}
(\phi,\Phi)\in \text{Mor}_{\Sf}[(\As_1,\H_1,D_1), (\As_2,\H_2,D_2)],  \\
 (\As_1,\H_1,D_1)\xrightarrow{(\phi,\Phi)} (\As_2,\H_2,D_2),
\end{gather*}
where $\phi: \As_1\to \As_2$ is a $*$-morphism between the pre-\cs s $\As_1,\As_2$ and $\Phi: \H_1\to \H_2$ is a bounded linear map in $\B(\H_1,\H_2)$ that ``intertwines'' the representations $\pi_1, \pi_2\circ\phi$ and the Dirac operators $D_1, D_2:$
\begin{gather}
\pi_2(\phi(x))\circ \Phi=\Phi\circ \pi_1(x), \quad \forall x \in \As_1, \label{eq: Phi}\\
D_2\circ\Phi= \Phi \circ D_1, \notag
\end{gather}
i.e.~such that the following diagrams commute for every $x\in \As_1:$
\begin{equation*}
\xymatrix{
{\H_1} \ar[d]_{D_1} \ar[r]^{\Phi}\ar@{}[dr]|{\circlearrowleft} & {\H_2}\ar[d]^{D_2} \\
{\H_1}\ar[r]^{\Phi} & {\H_2}
	}
\quad 
\xymatrix{
{\H_1} \ar[d]_{\pi_1(x)} \ar[r]^{\Phi}\ar@{}[dr]|{\circlearrowleft} & {\H_2}\ar[d]^{\pi_2\circ\phi(x)} \\
{\H_1}\ar[r]^{\Phi} & {\H_2}
	} 
\end{equation*}
Of course, the intertwining relation between the Dirac operators makes sense only on the domain of $D_1$. 
In the rest of the paper, we tacitly assume that $\Phi$ carries the domain of $D_1$ into that of $D_2$.

Note also that such a definition of morphism implies quite a strong relationship between the spectra of the Dirac operators of the two spectral triples. 

Loosely speaking, in the commutative case (see~\cite{BCL2} for details), one should expect such definition to become relevant only for maps that ``preserve the geodesic structures'' (totally geodesic immersions and totally geodesic submersions).
Clearly our definition of morphism contains as a special case the notion of (unitary) equivalence of spectral triples~\cite[pp.~485-486]{FGV}.

\subsection{Categories of Real and Even Spectral Triples.}

In the case of real spectral triples, we can define a natural notion of morphism simply by requiring that the morphisms be compatible with the real structures in the following sense: given two real spectral triples $(\As_j,\H_j,D_j,J_j),$ with $j=1,2,$ a morphism in our \emph{category of real spectral triples} $\Sf_r$ will be a morphism of spectral triples  
\begin{gather*}
(\phi,\Phi)\in \text{Mor}_{\Sf}[(\As_1,\H_1,D_1), (\As_2,\H_2,D_2)], \\
 (\As_1,\H_1,D_1)\xrightarrow{(\phi,\Phi)} (\As_2,\H_2,D_2),
\end{gather*}
such that $\Phi$ also ``intertwines'' the real structure operators $J_1, J_2:$
\begin{equation}\label{eq: a}
J_2\circ \Phi=\Phi \circ J_1,
\end{equation}
i.e.~such that the following diagram commutes:
\begin{equation*}
\xymatrix{
{\H_1} \ar[d]_{J_1} \ar[r]^{\Phi}\ar@{}[dr]|{\circlearrowleft} & {\H_2}\ar[d]^{J_2} \\
{\H_1}\ar[r]^{\Phi} & {\H_2}
	}
\end{equation*}

In a completely similar way, we can consider even spectral triples $(\As,\H,D,\Gamma)$ and 
define the \emph{category of even spectral triples} $\Sf_e,$ considering only those morphisms 
\begin{equation*}
(\As_1,\H_1,D_1)\xrightarrow{(\phi,\Phi)} (\As_2,\H_2,D_2),
\end{equation*}
such that $\Phi$ ``intertwines'' with the parity operators $\Gamma_1, \Gamma_2$, i.e. such that:
\begin{gather}\notag
\Gamma_2\circ \Phi=\Phi \circ \Gamma_1, \\ \label{eq: b}
\xymatrix{
{\H_1} \ar[d]_{\Gamma_1} \ar[r]^{\Phi}\ar@{}[dr]|{\circlearrowleft} & {\H_2}\ar[d]^{\Gamma_2} \\
{\H_1}\ar[r]^{\Phi} & {\H_2}
	}
\end{gather}

Again, in the case of real even spectral triples $(\As,\H,D,J,\Gamma)$ we will obtain a \emph{category of real even spectral triples} $\Sf_{re},$ choosing those morphisms that satisfy at the same time both the intertwining conditions~\ref{eq: a} and~\ref{eq: b} above.

Of course the category $\Sf_{re}$ of real even spectral triples is in general a non-full subcategory of both the categories $\Sf_r$ and $\Sf_e$ which are in turn non-full subcategories of $\Sf.$

\begin{remark}
According to our definition of morphisms, an automorphism of a real spectral triple $(\A, \H, D, J)$ in the categorical sense is given by a pair $(\phi, \Phi)$ with $\phi\in \Aut(\A)$ and $\Phi\in \B(\H)$ implementing $\phi$ and commuting with $D$ and $J$.
If we had required from the beginning the $\Phi$ appearing in~\eqref{eq: Phi} to be isometric, we would have obtained an extension of the isometry subgroup of $\text{Aut}^+$, 
the latter being the group of diffeomorphisms preserving the K-homology class of the spectral triple introduced by A.~Connes~\cite[Section~XI]{C5}. 
\end{remark}

Define $\Omega^p_D(\As):=\spa \{\pi(a_0)[D,\pi(a_1)]\cdots[D,\pi(a_p)] \ | \ a_0,\dots, a_p\in \As\},$ the space of $p$-forms. Every morphism 
$(\phi,\Phi): (\As_1,\H_1,D_1)\to(\As_2,\H_2,D_2)$
of spectral triples intertwines the $p$-forms according to the following formula:  
\begin{align*}
\Phi\circ &\sum_{j=1}^N \pi_1(a^{(j)}_0)[D_1,\pi_1(a^{(j)}_1)]\cdots[D_1,\pi_1(a^{(j)}_p)]= \\
&=\sum_{j=1}^N \pi_2(\phi(a^{(j)}_0))[D_2,\pi_2(\phi(a^{(j)}_1))]\cdots [D_2,\pi_2(\phi(a^{(j)}_p))]\circ \Phi. 
\end{align*}

\begin{remark}
Our morphisms of spectral triples are compatible with the inner deformation of the metric (see A.~Connes~\cite{C3,C4,C5}) in the following sense. 
Suppose that $(\phi,\Phi): (\As_1,\H_1,D_1)\to(\As_2,\H_2,D_2)$ is a morphism of spectral triples. Let us consider the two spectral triples $(\As_1, \H_1,D_1+A_1)$ and $(\As_2,\H_2,D_2+A_2)$ obtained by Morita ``self-equivalences'' of $\As_1$ and $\As_2$ using the ``gauge potentials'' $A_1\in \Omega^1_{D_1}(\As_1)$ and $A_2\in \Omega^1_{D_2}(\As_2)$,  respectively.
We notice that $(\phi,\Phi)$ continues to be a morphism of the ``deformed'' spectral triples 
if and only if $\Phi\circ A_1=A_2\circ \Phi$. 
\end{remark}

\section{Discrete Groups with Weights.} \label{sec: gr}

In order to prove the perfect mathematical naturality of our tentative definition of morphism of spectral triples, we provide here one interesting example of covariant functor with values in our category $\Sf.$

\subsection{Preliminaries on Group Algebras.}\label{sec: gr-alg}

For the benefit of the reader, we set up the framework by recalling a few properties of group algebras (of discrete groups) and their representations.

Let $G$ be a group equipped with the discrete topology\footnote{With this topology $G$ is of course a topological group.}. 

We recall that, given a group $G,$ we can always construct its \emph{group algebra} 
$\CC[G],$ that we will denote here by $\A_G:=\CC[G].$ \\ 
$\A_G$ consists of all the possible complex-valued functions on $G$ with finite support
$\A_G:=\{f: G \to \CC \ | \ f^{-1}\{\CC-\{0\}\} \ \text{is a finite set} \}$,
with sum and ``scalar'' multiplication by complex numbers defined pointwise: 
$(f+h)(x):= f(x)+g(x),$ 
$(\alpha f)(x):= \alpha(f(x)),$ 
and multiplication defined by the ``convolution'' product:
$(f*h)(z):=\sum_{\{(x,y)\ | \ xy=z\}} f(x)g(y).$
It is quickly established that $\A_G$, with the previously defined operations, is a complex associative unital algebra whose identity\footnote{Here $e$ denotes the identity element of $G$.} is $\delta^G_e(x)$, where $\delta^G_y(x):=\begin{cases} 1, \quad x=y \\0, \quad x\neq y \end{cases}$, and that $\A_G$ becomes a unital associative involutive algebra with the natural involution $(f^*)(x):=\cj{f(x^{-1})}$. 

\begin{proposition}\label{pr: g-a}
There exists a covariant functor $\A$ from the category $\Gf$ of groups with homomorphisms, to the category $\Af$ of associative complex unital involutive algebras with unital involutive algebra homomorphisms that to every group $G$ associates the group algebra $\A_G.$
\end{proposition}
\begin{proof} 
We have to define the functor on morphisms i.e.~given a homomorphism $\phi: H\to G$ between two groups $H$ and $G,$ we have to define a unital involutive homomorphism $\A_\phi:\A_G\to \A_H$ between  the group algebras. 

First of all notice that every group $G$ can be naturally ``embedded'' inside its group algebra by $\delta^G: G\to \A_G$, $z \mapsto \delta^G_z$.  
The map $\delta^G$ is injective, unital (i.e.,~$e\mapsto \delta^G_e$), multiplicative  (i.e.~$\delta^G_x*\delta^G_y=\delta^G_{xy}$), involutive (i.e.~$(\delta^G_x)^*=\delta^G_{x^{-1}}$). 

Then recall that for a given group $G,$ $(\A_G,\delta^G)$ is a free object over $G$ in the category of unital associative involutive algebras i.e.~every unital multiplicative involutive function $\psi: G\to \B$ from $G$ to a unital associative involutive algebra $\B$ can be ``lifted'' to a unital involutive algebra homomorphism $\Psi$ that makes the following diagram commutative:
\begin{equation*}
\xymatrix{G \ar[r]^{\delta^G} \ar[dr]_{\psi} \ar@{}[dr]^(0.5){\circlearrowleft}
& \A_G \ar[d]^{\Psi} \\
 & \B
	}
\end{equation*}
Finally take in the above diagram respectively $\B:=\A_H,$ $\psi: G\to \A_H$ defined by $\psi:= \delta^H\circ\phi$ in order to get the desired morphism of unital involutive algebras $\A_\phi:=\Psi.$ 
The association $\phi \mapsto \A_\phi$ is ``funtorial'' i.e.~respects compositions and identity functions. 
\end{proof}

\begin{proposition}\label{prop: hil}
On the complex vector space $\A_G$ there exists a natural inner product given by:
\begin{equation*}
\ip{f}{h}:=\sum_{x\in G} \cj{f(x)}h(x). 
\end{equation*}
With this inner product $\A_G$ is a pre-Hilbert space.
The completion of $\A_G$ with respect to the previous inner product is a Hilbert space.  
\end{proposition}

The Hilbert space constructed in proposition~\ref{prop: hil} is naturally identified with the Hilbert space $l^2(G):=L^2(G,\mu),$ where $\mu$ is the counting measure on the discrete group $G.$ In the following we will always denote this Hilbert space by $\H_G:= l^2(G).$

\begin{proposition}
There is a natural unital representation $\pi_G^0:\A_G\to \L(\A_G)$ of the group algebra $\A_G$ over itself by left action (by convolution). 
The representation is faithful.
\end{proposition}
\begin{proof} 

To every element $f\in \A_G$ we associate the element $\pi_G^0(f): \A_G\to \A_G$ given by
$(\pi_G^0(f))(h):=f*h,$ for every $h\in \A_G.$

From the definition of $\pi_G^0$ it is clear that $\pi_G^0(f)\in \L(\A_G)$ and that $f\mapsto \pi_G^0(f)$ is a linear function: $\pi_G^0 \in \L(\A_G; \L(\A_G)).$

By direct calculation, $\pi_G^0$ is multiplicative and unital hence a representation.

The injectivity of $\pi_G^0$ follows from the triviality of the kernel
(as in any unital left-regular representation): if $f\neq 0,$ then $\pi_G^0(f)(\delta_e^G)=f*\delta^G_e=f\neq 0.$
\end{proof}

\begin{corollary}\label{cor: rep}
There is a natural faithful representation $\pi_G: \A_G \to \B(\H_G)$ of the group algebra $\A_G$ as bounded operators on the Hilbert space $\H_G.$ 
\end{corollary}
\begin{proof} 

The operator $\pi_G^0(f)\in \L(\A_G)$ is a bounded operator on the pre-Hilbert space $\A_G.$
To prove this note that if $f=\sum_{x\in G}f(x)\delta_x^G,$ by the linearity of $\pi_G^0,$ we have $\pi_G^0(f)=\sum_{x\in G}f(x)\pi_G^0(\delta_x^G)$ so that it is enough to prove the boundedness of the operators $\pi_G^0(\delta_x^G)$ for all $x\in G.$ This follows immediately from the fact that $\pi_G^0(\delta_x^G)$ is an isometry of the inner product space $\A_G:$ 
\begin{equation*}
\|\pi_G^0(\delta_x^G)(h)\|^2=\|\delta_x^G*h\|^2=\sum_{z\in G}|h(x^{-1}z)|^2=\sum_{z\in G}|h(z)|^2= \|h\|^2.
\end{equation*}
By linear extension theorem, $\pi_G^0(f)$ extends to a bounded operator $\pi_G(f)$ on $\H_G$ with the same norm. 
\end{proof}

The representation $\pi_G: \A_G\to\B(\H_G),$ in
corollary~\ref{cor: rep}, is nothing but the \emph{left-regular representation} $\lambda_G: \CC[G]\to \B(l^2(G)).$

\begin{proposition}\label{pr: j}
The exists a natural antilinear involution $J_G:\H_G\to \H_G$ on the Hilbert space $\H_G.$
\end{proposition}
\begin{proof}
On the pre-Hilbert space $\A_G,$ the algebra involution $*: \A_G\to \A_G,$ defined by $f^*(x):=\cj{f(x^{-1})},$ is antilinear and isometric:
\begin{equation*}
\ip{f^*}{g^*}=\sum_{x\in G}\cj{\cj{f(x^{-1})}}\ \cj{g(x^{-1})}=\sum_{x\in G}\cj{g(x^{-1})}f(x^{-1})=\sum_{x\in G}\cj{g(x)}f(x)=\ip{g}{f}.
\end{equation*}
By linear extension theorem (for antilinear maps), there exists a unique antilinear extension $J_G:\H_G\to \H_G$ to the closure $\H_G$ of $\A_G.$ 
The map $J_G$ is antilinear, involutive, isometric.
\end{proof}

\subsection{Preliminaries on Weighted Groups.}

\begin{definition}
By a \emph{weight} on a group $G$ we mean a real-valued function $\omega: G \to \RR.$ 
Given two weighted groups $(G,\omega_G)$ and $(H,\omega_H)$, we say that a function
$\phi: G\to H$ is a \emph{weighted homomorphism} if:
\begin{gather*}
\phi: G \to H \quad \text{is a group homomorphism and}\quad \omega_G=\phi^\bullet(\omega_H):=\omega_H \circ \phi.
\end{gather*} 
A weight is called \emph{proper}
if for every $k\in \NN,$ $\omega_G^{-1}([-k,+k])$ is a finite set in $G$. 
\end{definition}

Note that proper weights exist only on countable groups.

\begin{remark}
A special case of weight on a group $G,$ is given by the notion of a \emph{length function} on a group~\footnote{Here we follow the definition used by M.~Rieffel~\cite[Section~2]{Ri1}.}~\cite{C1} i.e.~a function $\ell_G: G\to \RR$ such that:
\begin{gather*}
\ell_G(xy)\leq \ell_G(x)+\ell_G(y), \quad \forall x,y \in G, \\
\ell_G(x^{-1})=\ell_G(x), \quad \forall x\in G, \\
\ell_G(x)=0 \iff x=e, \quad \text{where $e\in G$ is the identity element of $G$}.
\end{gather*}
Of course a length function is always positive since: \\ $0=\ell_G(e)=\ell_G(xx^{-1})\leq\ell_G(x)+\ell_G(x^{-1})=2\ell_G(x)$ 
for all $x\in G.$  

A weighted homomorphism of groups with length is called an \emph{isometry}.
The previous conditions actually imply that every isometry is injective:
\begin{equation*}
\phi(x)=e_H \imp \ell_H(\phi(x))=\ell_H(e_H)=0 \imp \ell_G(x)=0 \imp x=e_G.
\end{equation*}
\end{remark}

\begin{proposition}
The class of (proper) weighted groups with weighted homomorphisms forms a category. 
The class of groups equipped with a (proper) length function when the morphisms are the isometries, is a full subcategory. 
\end{proposition}
\begin{proof}
The composition of weighted homomorphisms (respectively isometries) $\phi: G\to H$ and $\psi: H\to K$ is a weighted homomorphism (isometry):
\begin{equation*}
(\psi\circ\phi)^\bullet(\omega_K)=\phi^\bullet(\psi^\bullet(\omega_K))=\phi^\bullet(\omega_H)=\omega_G.
\end{equation*}
For every object $(H,\omega_H),$ the identity isomorphism $\iota: H\to H$ is a weighted homomorphism (isometry) that satisfies $\psi\circ\iota = \psi,$ and $\iota\circ \phi=\phi$ for every composable weighted homomorphisms $\phi, \psi$.
\end{proof}

Of course the category of normed spaces with linear norm-preserving maps is a (non-full) subcategory of the category of abelian groups with length function (the length function being the norm) and isometries as defined above coincide with the well-known concept of norm-preserving maps in normed spaces. 

\begin{proposition}\label{prop: a-h}
There is a covariant functor $\A$ from the category $\Gf_i$ of groups with injective homomorphism as arrows, to the category of $p\Hf_i$ pre-Hilbert spaces with isometries.

In the same way, we have a covariant functor $\H$ from the category $\Gf_i$ of groups with injective morphism to the category $\Hf_i$ of Hilbert spaces with isometries.

The functors $\A$ and $\H$ are left exact.
\end{proposition}
\begin{proof} 
The functor on objects is defined by $G\mapsto \A_G\in p\Hf_i$ and by $G\mapsto \H_G\in \Hf_i$ respectively.

To define the functor on morphisms, we first note that for any given group $G,$ the set $\{\delta^G_x \mid x\in G\}$ is a (Hamel) basis for the vector space $\A_G$ that is orthomormal with respect to the inner product in $\A_G.$

If the function $\phi: G \to H$ is a monomorphism, the induced (linear) map 
$\A_\phi: \A_G \to \A_H$ becomes an isometry because it maps $\delta^G_x$ to $\delta^H_{\phi(x)}$ i.e.~it sends an orthonormal basis to an orthonormal set.

Since $\A_\phi$ is an isometry, it is bounded as a map from $\A_G$ to $\H_H$ and it can be uniquely extended to an isometry $\H_\phi: \H_G \to \H_H.$

The associations $\phi \mapsto \A_\phi$ and $\phi \mapsto \H_\phi$ satisfy all the functorial properties.
\end{proof}

The following theorem is a well-known result of A.~Connes~\cite[Lemma~5]{C1}:
\begin{theorem}\label{th: co-sptr}
To every pair $(G,\omega_G)$ where $G$ is a discrete countable group and $\omega_G$ is a 
weight function on $G,$ we can associate a triple $(\A_G, \H_G, D_{\omega_G})$ given as follows:
\begin{itemize}
\item
$\A_G$ is the group algebra of $G$ as defined above in subsection~\ref{sec: gr-alg}.
\item
$\H_G$ is the Hilbert space of $G$ as defined above in proposition~\ref{prop: hil}.
\item
The representation of the algebra $\A_G$ on $\H_G$ is the left-regular representation $\pi_G: \A_G\to \B(\H_G)$ defined above in corollary~\ref{cor: rep}.
\item
The Dirac operator $D_{\omega_G}$ is the pointwise multiplication operator by the 
weight function $\omega_G$, i.e.
\begin{equation*}
(D_{\omega_G} \xi)(x):= \omega_G(x)\xi(x), \quad \forall x\in G ,
\end{equation*} 
naturally defined on the domain $\{\xi \in \H_G \ | \ \sum_{x \in G} |\omega_G(x) \, \xi(x)|^2 < \infty \}$.
\end{itemize}

The triple $(\A_G,\H_G,D_{\omega_G})$ is a spectral triple if and only if the weight $\omega_G$ is proper and such that, for all $x\in G,$ the differences\footnote{Where $\tau_x(\omega_G): y\mapsto \omega_G(x^{-1}y)$ is the ``left $x$-translated'' of $\omega_G.$} $[\omega_G -\tau_x(\omega_G)]: G\to \RR,$ are bounded real-valued functions.
\end{theorem}
\begin{proof}
$\A_G$ is a pre-C$^*$ algebra: defining $\|f\|:=\|\pi_G(f)\|_{\H_G},$ we see that the C$^*$-property $\|f^*f\|=\|f\|^2$ is immediate\footnote{It must be pointed out that, denoting by $\text{C}^*_r(G)$ the closure of $\A_G$ in the norm defined above, the correspondence $G\mapsto \text{C}^*_r(G)$ is not functorial, in general.
It becomes so, for the full subcategory of amenable groups. 
In the case of non amenable groups we do not have finite dimensional spectral triples 
(see A.~Connes~\cite[Theorem~19]{C1}).}.

The Dirac operator $D_{\omega_G}$ is self-adjoint and has compact resolvent if and only if $\omega_G$ is proper. 

Every element $f\in \A_G$ can be written as $f= \sum_{x\in G} f(x)\delta^G_x.$

It follows that $\pi_G(f)= \sum_{x\in G}f(x)\pi_G(\delta^G_x)$ and we have:
\begin{equation*}
\|[D_{\omega_G}, \pi_G(f)]\|=\|\sum_{x\in G}f(x)[D_{\omega_G}, \pi_G(\delta^G_x)]\| \leq \sum_{x\in G}|f(x)|\cdot\| [D_{\omega_G}, \pi_G(\delta^G_x)]\|,
\end{equation*}
so that, in order to show the boundedness of $[D_{\omega_G},\pi_G(f)]$ it is enough to show the boundedness of $[D_{\omega_G},\pi_G(\delta^G_x)]$ for all $x\in G.$ 

Now, from the fact that $\pi_G(\delta^G_x)$ is unitary in $\H_G,$ we have:
\begin{align*}
\|[D_{\omega_G},\pi_G(\delta^G_x)]\|&=\|D_{\omega_G}\pi_G(\delta^G_x)-\pi_G(\delta^G_x)D_{\omega_G}\| \\
&=\|(D_{\omega_G}-\pi_G(\delta^G_x)D_{\omega_G}\pi_G(\delta^G_x)^{-1})\pi_G(\delta^G_x)\| \\
&=\|D_{\omega_G}-\pi_G(\delta^G_x)D_{\omega_G}\pi_G(\delta^G_x)^{-1}\|
\end{align*}
and since, by direct calculation, we get $\pi_G(\delta^G_x)D_{\omega_G}\pi_G(\delta^G_x)^{-1}=D_{\tau_x(\omega_G)},$ 
where $\tau_x(\omega_G): y\mapsto \omega_G(x^{-1}y),$ we see that $\|[D_{\omega_G},\pi_G(\delta^G_x)]\|=\|D_{\omega_G}-D_{\tau_x(\omega_G)}\|.$
Since $\|D_{\omega_G}-D_{\tau_x(\omega_G)}\|=\|\omega_G-\tau_x(\omega_G)\|_\infty:=\sup\{|\omega_G(y)-\omega_G(x^{-1}y)| :  y\in G\},$ the assertion is proved\footnote{Let $u: G\times G \to \TT$ be a normalized $2$-cocycle on $G$ and consider the
left-regular representation of $G$ twisted by $u$, defined by
$(\pi^G_u(\delta^G_x))\delta^G_y:=u(y^{-1}x^{-1}, x)\delta^G_{xy}$. Then up to minor modifications the same argument shows that the Dirac operator $D_{\omega_G}$ also gives a spectral triple over the ``twisted group algebra'' generated by the $\pi^G_u(\delta_x^G)$' s.}.
\end{proof}

\begin{remark}
In the case of length funtions on groups, the last condition of theorem~\ref{th: co-sptr} is automatically satisfied:
\begin{equation*}
\ell_G(y)-\tau_x(\ell_G)(y)=\ell_G(y)-\ell_G(x^{-1}y)\leq\ell_G(x).
\end{equation*}
\end{remark}

\begin{lemma}\label{lem: core}
$\A_G\subset \H_G$ is an invariant core for the operator $D_{\omega_G}.$
\end{lemma}
\begin{proof}
Suppose that\footnote{Operators and their graphs are denoted with the same symbol.} 
$(\xi,\eta)\in D_{\omega_G}.$ Since $\A_G$ is dense in $\H_G,$ there is a sequence $\xi_n\xrightarrow{n\to \infty}\xi$ with $\xi_n \in \A_G.$
We show that it is possible to choose the sequence $\xi_n\in \A_G$ in such a way that $D_{\omega_G}(\xi_n)\xrightarrow{n\to \infty} \eta$.
In fact, selecting an arbitrary well ordering $n\mapsto x_n\in G$ in the support set of $\xi,$ we can always define $\xi_n:=\sum^n_{j=0}\xi(x_j)\delta^G_{x_j}$ and check that $\xi_n\xrightarrow{n\to \infty}\xi$ and also $D_{\omega_G}(\xi_n)\xrightarrow{n\to \infty}\eta$ so that $(\xi,\eta)\in \cj{D_{\omega_G}|_{\A_G}}$ i.e.~$\A_G$ is a core for $D_{\omega_G}.$

Of course, since $\xi_n$ has finite support, $D_{\omega_G}(\xi_n)$ also has finite support and so $D_{\omega_G}(\xi_n)\in \A_G$. In particular $\A_G$ is an invariant subspace for $D_{\omega_G}.$
\end{proof}

\begin{lemma}\label{lem: l-t}
Given the weight $\omega_G: G\to \RR$ on the group $G,$ the following conditions are equivalent\footnote{By definition, $\tau'_x(\omega_G): y\mapsto \omega_G(yx^{-1})$ is the ``right'' translation of $\omega_G$ by $x.$}:
\begin{gather*}
\forall x\in G \quad \omega_G -\tau_x(\omega_G) \quad \text{is constant}; \\
\omega_G=\alpha+\phi, \quad \text{where $\alpha$ is a constant and $\phi:G\to \RR$ is a homomorphism;}\\
\forall x\in G, \quad \omega_G - \tau'_x(\omega_G) \quad \text{is constant}; \\
\omega_G(xzy^{-1})-\omega_G(zy^{-1})=\omega_G(xz)-\omega_G(z), \quad \forall x,y,z\in G.
\end{gather*}
\end{lemma}
\begin{proof}
By direct calculation if $\omega_G=\alpha+\phi$ then $\omega_G-\tau_x(\omega_G)$ and $\omega_G-\tau_x'(\omega_G)$ are constant. That $\omega_G-\tau_x(\omega_G)$ being constant is equivalent to 
\begin{equation}\label{eq: 1}
\omega_G(xg)=\omega_G(g)-\phi(x^{-1}),
\end{equation}
for some function $\phi: G\to \RR$. Taking $x=g^{-1}$ in the previous equation we have $\phi(g)=\omega_G(g)-\omega_G(e_G)$. Hence equation~\eqref{eq: 1} implies 
\begin{equation} \label{eq: 1.1}
   \phi(xg) = \phi(g)-\phi(x^{-1}) 
\end{equation}
and (taking $g=e_G$) $\phi(x)=-\phi(x^{-1})$. Substituting this in equation~\eqref{eq: 1.1}, we see that $\phi$ is a homomorphism so that $\omega=\alpha+\phi$ with $\alpha:=\omega(e_G)$.  
The same proof applies to the case $\omega_G-\tau_x'(\omega_G)$ being constant.

The last equation is easily reduced to equivalence to the constancy of $\omega_G-\tau_x'(\omega_G)$ by substitutions.
\end{proof}

In view of their relevance for the construction of spectral triples, weights satisfying the last condition in theorem~\ref{th: co-sptr} deserve a special name.
\begin{definition}
A weight $\omega_G$ on the group $G$ is said to be a \emph{Dirac weight} if  
$\omega_G - \tau_x(\omega_G)$ are bounded functions, for every $x \in G$.
\end{definition}

The following proposition is essentially a restatement of the results already obtained by M.~Rieffel~\cite[See the end of Section~2]{Ri2}.

\begin{proposition}
Given a proper Dirac-weighted countable group $(G,\omega_G),$ the anti-unitary operator $J_G$ defined in proposition~\ref{pr: j} is a real structure on the spectral triple $(\A_G, \H_G,D_{\omega_G})$ if and only if either $\omega_G$ is a constant function or $\omega_G$ is a homomorphism of groups.
\end{proposition}
\begin{proof}
We have $J_G^2=\id_{\H_G}.$ 

By linear extension, the condition 
$(J_G D_{\omega_G})(\xi)=\pm (D_{\omega_G}J_G)(\xi)$ for $\xi\in\H_G$ holds if and only if $(J_G D_{\omega_G})(\delta^G_x)=\pm (D_{\omega_G}J_G)(\delta^G_x)$, which is also equivalent to:
\begin{equation}\label{eq: pm}
\omega_G(x^{-1})=\pm \omega_G(x)\quad \forall x\in G.
\end{equation}
There is no problem at all to verify the property
\begin{equation}\label{eq: r1}
(J_G\pi_G(g)J_G)\circ \pi_G(f)(\xi)=\pi_G(f)\circ (J_G\pi_G(g)J_G)(\xi), \quad \forall \xi\in \H_G.
\end{equation}
In fact, for all $f,g,\xi\in \A_G:$
\begin{align*}
[\pi_G(f)\circ &(J_G\pi_G(g)J_G)](\xi)=f*J_G(g*(J_G(\xi)))=f*(J_G^2(\xi)*J_G(g))=\\
				&=f*\xi*(J_G(g))=J_G(g*(J_G(\xi))*(J_G(f)))=\\
				&=(J_G\pi_G(g)J_G)(f*\xi)=[(J_G\pi_G(g)J_G)\circ\pi_G(f)](\xi)
\end{align*}
and by linear extension theorem (since $J_G$, $\pi_G(f)$, and $\pi_G(g)$ are all bounded) condition~\eqref{eq: r1} holds for all $\xi\in \H_G.$

We now prove that the first order condition 
\begin{equation}\label{eq: r2}
[D_{\omega_G},\pi_G(f)]_- \circ  (J_G\pi_G(g)J_G)(\xi)
=(J_G\pi_G(g)J_G)\circ [D_{\omega_G},\pi_G(f)]_-(\xi),  
\end{equation}
for all $f,g\in \A_G$ and all $\xi\in \H_G,$ holds if and only if $\omega_G-\tau_x(\omega_G)$ are constant functions.
Since all the operators involved are bounded on the Hilbert space $\H_G,$ by linear extension theorem, it is enough to check the first order condition only for every $\xi\in \A_G.$

Let $f=\sum_{x\in G}f(x)\delta^G_x,$ $g=\sum_{y\in G}g(y)\delta^G_y$ and $\xi=\sum_{z\in G}\xi(z)\delta^G_z$ be three elements in $\A_G.$ Substitution in equation~\eqref{eq: r2} above and (anti-)linearity yield
\begin{align*}
&\sum_{x,y,z\in G}\cj{f(x)}g(y)\xi(z)\cdot [D_{\omega_G}, \delta_x^G]_-\circ(J_G\pi_G(\delta^G_y) J_G) (\delta_z^G)
\\
&=\sum_{x,y,z\in G}\cj{f(x)}g(y)\xi(z) \cdot (J_G\pi_G(\delta^G_y) J_G)\circ [D_{\omega_G}, \delta_x^G]_-(\delta_z^G).
\end{align*}
This last equation holds if and only if, for all $x,y,z\in G:$
\begin{equation*}
[D_{\omega_G}, \delta_x^G]_-\circ(J_G\pi_G(\delta^G_y) J_G) (\delta_z^G)=
(J_G\pi_G(\delta^G_y) J_G)\circ [D_{\omega_G}, \delta_x^G]_-(\delta_z^G).
\end{equation*}
By direct calculation we have:
\begin{gather*}
[D_{\omega_G}, \delta_x^G]_-\circ(J_G\pi_G(\delta_y^G) J_G) (\delta_z^G)= 
\omega_G(xzy^{-1})-\omega_G(zy^{-1})\delta^G_{xzy^{-1}},\\
(J_G\pi_G(\delta_y^G) J_G)\circ [D_{\omega_G}, \delta_x^G]_-(\delta_z^G)=
\omega_G(xz)-\omega_G(z)\delta^G_{xzy^{-1}}.
\end{gather*}
Hence our result is that the first order condition~\eqref{eq: r2} holds if and only if 
\begin{equation*}
\omega_G(xzy^{-1})-\omega_G(zy^{-1})=\omega_G(xz)-\omega_G(z), \quad \forall x,y,z\in G.
\end{equation*}
and this, by lemma~\ref{lem: l-t}, is equivalent to the fact that $\omega_G=\alpha+\phi$ 
where $\alpha: G\to \RR$ is constant and $\phi: G\to\RR$ is a homomorphism of groups.

Now, equation~\eqref{eq: pm} above, in the plus case, is equivalent to $\phi=0$ and so to  $\omega_G=\alpha$ being a constant. In the minus case, it is equivalent to $\alpha=0$ and so to  $\omega_G=\phi$ being a homomorphism of groups.
\end{proof}

\begin{remark} 
The spectral triple $(\A_G,\H_G,D_{\omega_G})$ is regular (see M.~Rieffel~\cite[End of section~2]{Ri2}).
For instance, $[|D_{\omega_G}|,\pi_G(f)]$ and $[|D_{\omega_G}|, [D,\pi_G(f)]\,]$ are bounded for all $f\in \A_G$ 
as a consequence of the following estimates which can be obtained by repeating the argument in the proof of theorem~\ref{th: co-sptr}:
\begin{gather*}
\| \, [D_{|\omega_G|},\pi_G(\delta^G_x)] \, \| \leq
\|\tau_x(|\omega_G|)-|\omega_G|\, \|_\infty \leq \|\tau_x(\omega_G)-\omega_G\|_\infty; \\
\| \, [\, D_{|\omega_G|}, [\, D_{\omega_G},\pi_G(\delta^G_x)]\,] \, \| \leq  \|\tau_x(\omega_G)-\omega_G\|_\infty^2,
\end{gather*}
and more generally
\begin{gather*}
\| \, [\, |D_{\omega_G}|,\cdots, [\, |D_{\omega_G}|,\pi_G(\delta^G_x)
\underbrace{]\,\cdots ]}_n \, \| \leq  \|\tau_x(\omega_G)-\omega_G\|_\infty^n.\\
\| \, [\, |D_{\omega_G}|,\cdots, [\, |D_{\omega_G}|,[D_{\omega_G},\pi_G(\delta^G_x)]\,\underbrace{]\,\cdots ]}_n \, \| \leq  \|\tau_x(\omega_G)-\omega_G\|_\infty^{n+1}.\\
\end{gather*}
\end{remark}

\begin{remark}\label{rk: gr-sp-tr}

On the real spectral triple $(\A_G,\H_G,D_{\omega_G},J_G),$ it is impossible to introduce a grading operator $\Gamma_G,$ (unless $\omega_G$ is the zero function\footnote{In the case $\omega_G$ equal to zero, a convenient grading is given by $\Gamma(\delta^G_x):=\delta^G_{x^{-1}}.$}).
This is because if $\omega_G$ is a non-zero constant the equation $D_{\omega_G}\Gamma=-\Gamma D_{\omega_G}$ cannot be satisfied. On the other hand, if $\omega_G$ is a homomorphism, then we are in the case $J_G D_{\omega_G}=-D_{\omega_G}J_G$ which, from the table at the end of section~\ref{tb: J}, is incompatible with the existence of a grading. Of course, ``doubling'' in an appropriate way the Hilbert space $\H_G,$ we can easily get another graded real spectral triple: 
\begin{itemize}
\item
the pre-C$^*$-algebra is the same group algebra $\A_G;$
\item
the Hilbert space is given by the direct sum $\H_G\oplus \H_G;$
\item
the representation of $\A_G$ in $\H_G\oplus \H_G$ is the direct sum representation $\pi_G\oplus \pi_G$ i.e.~for all $f\in \A_G$ and $\xi,\eta\in \H_G:$
\begin{equation*}
[\pi_G\oplus\pi_G(f)](\xi\oplus \eta):=\begin{bmatrix}
											\pi_G(f)	& 0 \\
											0 & \pi_G(f)
											\end{bmatrix}
											\cdot
											\begin{bmatrix}
											\xi \\
											\eta
											\end{bmatrix}
											=
											\begin{bmatrix}
											[\pi_G(f)](\xi) \\
											[\pi_G(f)](\eta)
											\end{bmatrix};
\end{equation*}
\item
the Dirac operator is given by:
\begin{equation*}
D_{\omega_G}\oplus (-D_{\omega_G})=
\begin{bmatrix}
D_{\omega_G} & 0 \\
0 & -D_{\omega_G}
\end{bmatrix};
\end{equation*}
\item
the grading operator is given by:
\begin{equation*}
\Gamma_G:=
\begin{bmatrix}
0 & 1 \\
1& 0
\end{bmatrix};
\end{equation*}
\item
the real structure is given by:
\begin{equation*}
J_G\oplus J_G =
\begin{bmatrix}
J_G & 0 \\
0 & J_G
\end{bmatrix}.
\end{equation*}
\end{itemize}
\end{remark}

\subsection{The Functor: Monomorphism Case.}\label{sec: fun}

\begin{theorem}\label{th: im}
There exists a covariant functor, from the category $\Gf^\omega_i$ of proper Dirac-weighted countable groups with weighted monomorphisms to the category of spectral triples, that to every $(G,\omega_G)$ associates $(\A_G,\H_G,D_{\omega_G}).$
\end{theorem}
\begin{proof}
We only need to prove existence of a functor on monomorphisms $\phi: G\to H.$ 
It is our purpose to show that the pair $(\A_\phi,\H_\phi)$ defined in proposition~\ref{prop: a-h} is a morphism 
\begin{equation*}
(A_G, \H_G, D_{\omega_G})\xrightarrow{(\A_\phi,\H_\phi)}(A_H, \H_H, D_{\omega_H})
\end{equation*}
of spectral triples.

This amounts to showing that for every $f\in \A_G$ and for every $\xi\in \H_G:$
\begin{equation*}
\H_\phi \circ \pi_G(f) (\xi)= \pi_H(\A_\phi(f))\circ \H_\phi (\xi);
\end{equation*}
and that:
\begin{equation}\label{eq: d-c}
\H_\phi\circ D_{\omega_G} (\xi)=D_{\omega_H}\circ \H_\phi(\xi).
\end{equation}

The first property follows from the fact that, for every $f\in \A_G$ and for every $\xi\in \A_G\subset \H_G$ we have:
\begin{equation*}
(\H_\phi\circ\pi_G(f))(\xi)=\H_\phi(f*\xi)=\A_\phi(f*\xi)=\A_\phi(f)*\A_\phi(\xi)
=\pi_H(\A_\phi(f))\circ \H_\phi(\xi).
\end{equation*}
Since, for every $f\in \A_G,$ the bounded operators $\H_\phi\circ\pi_G(f)$ and $\pi_H(\A_\phi(f))\circ \H_\phi$ coincide on the dense subspace $\A_G$ of $\H_G,$ the identity follows.

The second property is obtained from the fact that, for every $\xi\in \A_G\subset \H_G:$
\begin{align*}
\H_\phi\circ D_{\omega_G}(\xi) 
&= \H_\phi\left(\sum_{z\in G}\omega_G(z)\xi(z)\delta^G_z\right) 
								=\sum_{z\in G}\omega_G(z)\xi(z)\H_\phi(\delta^G_z) \\
&=\sum_{z\in G}\omega_G(z)\xi(z)\delta^H_{\phi(z)}
								= \sum_{z\in G}\omega_H(\phi(z))\xi(z)\delta^H_{\phi(z)} \\
&=D_{\omega_H}\left(\sum_{z\in G}\xi(z)\delta^H_{\phi(z)}\right) 
								=D_{\omega_H}\left(\sum_{z\in G}\xi(z)\H_\phi(\delta^G_z)\right) \\
&=D_{\omega_H}\circ \H_\phi\left(\sum_{z\in G}\xi(z)\delta^G_z\right) 
								=D_{\omega_H}\circ \H_\phi(\xi),
\end{align*}
so that the two operators $\H_\phi\circ D_{\omega_G}$ and $D_{\omega_H}\circ \H_\phi$ coincide on the dense subspace $\A_G$ of $\H_G.$
From the fact that $\A_G$ is an invariant subspace for $D_{\omega_G}$ and $\H_\phi,$ and from lemma~\ref{lem: core} above, we see that $\A_G$ is a core for both operators and the equality~\eqref{eq: d-c} follows.
\end{proof}

\begin{proposition}
Under the same assumptions as in theorem~\ref{th: im}, if the weight $\omega_G$ is a group homomorphism or a constant, then $(\A_\phi,\H_\phi)$ is a morphism of real spectral triples i.e.:
\begin{equation*}
\H_\phi \circ J_G=J_H\circ \H_\phi.
\end{equation*}
\end{proposition}
\begin{proof}
For every element $\xi=\sum_{x\in G}\xi(x)\delta_x^G\in \A_G\subset \H_G$ we have: 
\begin{align*}
\H_\phi\circ J_G\left( \sum_{x\in G} \xi(x)\delta_x^G \right) 
&= \sum_{x\in G} \H_\phi\circ J_G (\xi(x)\delta^G_x) = \sum_{x\in G} 			
				\cj{\xi(x)}\H_\phi(J_G(\delta^G_x))= \\
&= \sum_{x\in G}\cj{\xi(x)}\H_\phi(\delta^G_{x^{-1}}) = 
				\sum_{x\in G}\cj{\xi(x)}\delta^H_{\phi(x)^{-1}}= \\
&= \sum_{x\in G} J_H(\xi(x)\delta^H_{\phi(x)})=\sum_{x\in G} J_H\circ \H_\phi(\xi(x)\delta^G_x)= \\
&= J_H\circ \H_\phi\left( \sum_{x\in G} \xi(x)\delta^G_x\right).
\end{align*}
This means that both the operators $\H_\phi \circ J_G$ and $J_H\circ \H_\phi$ coincide on the dense subspace $\A_G$ of the Hilbert space $\H_G$ and, since they are bounded, it follows by linear extension theorem, that they are equal on all of $\H_G.$
\end{proof}

\begin{remark}
With the same notations developed in remark~\ref{rk: gr-sp-tr},
it is easily established that $(\A_\phi, \H_\phi\oplus\H_\phi)$ is a morphism of graded spectral triples (with real structure, when available). 
The association $G\mapsto \A_G,$ $\phi\mapsto (\A_\phi, \H_\phi\oplus\H_\phi)$ is functorial from the category $\Gf^\omega_i$ to the category $\Sf$ of spectral triples.
\end{remark}

An automorphism $\alpha$ of $G$ induces by functoriality an automorphism $\A_\alpha$ of the group algebra $\A_G$ implemented by the unitary $\H_\alpha$ on the Hilbert space $\H_G$
and, if $\alpha$ is also weighted, $(\A_\alpha,\H_\alpha)$ is an automorphism of the  spectral triple $(\A_G,\H_G,D_{\omega_G})$.
In particular, if $\alpha:=\ad_g,\ g\in G,$ is inner, $\H_\alpha = \pi_G(g)J\pi_G(g)J$
is the image of $g$ through the inner regular representation of $G$.\footnote{Note that $(\A_G,\H_G,D_{\omega_G})$, with $\A_G$ acting on $\H_G$ by the (linearization of the) inner regular representation of $G$, is a spectral triple too.}

\medskip

Equivalence classes of monomorphisms categorically correspond to subobjects, in our case,  subgroups.
Every subgroup $H$ of the weighted group $(G,\omega_G)$ comes naturally equipped with a weight function $\omega_H:=\omega_G|_H$ obtained by restriction of the original weight function on $G$ and the inclusion map $\iota: H\to G$ is a morphism in $\Gf^\omega_i$.
By proposition~\ref{prop: a-h} and theorem~\ref{th: im}, $(\A_{\iota},\H_\iota)$ is a monomorphism from the spectral triple $(A_H,\H_H,D_{\omega_H})$ to the spectral triple $(A_G,\H_G,D_{\omega_G})$. Similarly, one has the following: 
\begin{corollary}
The functor $\Fg:\Gf^\omega_i\to\Sf$ is left exact: every monomorphism of groups gives rise to a monomorphism of spectral triples. 
\end{corollary}

Note however that the functor $\Fg$ is not full: there are morphisms (even monomorphisms) of spectral triples over group algebras that are not obtained from monomorphisms of groups.
This fact might call for suitable modifications of our setting that could entail better functorial correspondences.

\subsection{Preliminaries on Charged Groups and Co-isometries.}

Before proceeding further, we need to collect a few more facts about weights and lengths on groups.
\begin{definition}
A \emph{charged group} is a weighted group $(G,\omega_G)$ such that the function 
$|\omega_G|:x\mapsto |\omega_G(x)|$ is a length function on $G.$

A homomorphism $\phi: G\to H$ between charged groups is called \emph{isometric} if $\phi^\bullet(|\omega_H|)=|\omega_G|.$
\end{definition}

\begin{remark}
Every group with length function $(G,\ell_G)$ is a charged group.

An isometric homomorphism between charged groups is continuous with respect to the metric topologies induced by the length functions. 

Every weighted homomorphism $\phi: G\to H$ between charged groups $(G,\omega_G)$ and $(H,\omega_H)$ is  isometric\footnote{Of course, $\phi$ is continuous.}.

The category of charged groups with isometric weighted morphisms is a full subcategory of the category of weighted groups with weighted monomorphisms. 
\end{remark}

\begin{definition}
Let $(G,\omega_G)$ and $(H,\omega_H)$ be two weighted groups.  A homomorphism of groups 
$\phi: G\to H$ is called a \emph{co-weighted homomorphism} if there exists a weighted homomorphism 
$\rho: H \to G$ such that $\phi\circ\rho=\iota_H.$ A co-weighted homomorphism between two charged groups 
is said to be \emph{co-isometric} if 
$|\omega_H(\phi(g))|\leq |\omega_G(g)|$ for all $g\in G$. 
\end{definition}

\begin{lemma}
Let $(G,\ell_G)$ be a group with length function and let $H$ be a normal subgroup of $G.$
The function $\ell_{G/H}: G/H \to \RR$ defined by
\begin{equation*}
\ell_{G/H}(xH):=\inf \{\ell_G(xh)\ | \ h\in H\}
\end{equation*}
is a length on $G/H$ called the \emph{quotient length.}
\end{lemma}
\begin{proof}

Using the normality of $H$ in $G$ we see that $\{xyhk\ | \ h,k\in H\}=\{xhyk\ | \ h,k\in H\}$. Hence this calculation follows:
\begin{gather*}
\ell_{G/H}(xyH)\leq \ell_G(xhyk)\leq \ell_G(xh)+\ell_G(yk) \quad \forall h,k\in H \imp \\
\ell_{G/H}(xyH) - \ell_G(xh)\leq \ell_G(yk), \quad \forall h,k\in H \imp \\
\ell_{G/H}(xyH) - \ell_G(xh)\leq \ell_{G/H}(yH), \forall h\in H \imp \\
\ell_{G/H}(xyH) -\ell_{G/H}(yH)\leq  \ell_G(xh), \quad \forall h\in H \imp \\
\ell_{G/H}(xyH) -\ell_{G/H}(yH)\leq  \ell_{G/H}(xH).
\end{gather*}

Since $H$ is normal in $G,$ we have $\{x^{-1}h \ | \ h\in H\}=\{(xh)^{-1}\ | \ h\in H\}$
and so:
\begin{align*}
\ell_{G/H}(x^{-1}H)&=\inf\{\ell_G(x^{-1}h)\mid h\in H\}=\inf\{\ell_G((xh)^{-1}) \ | \ h \in H\} \\ 
&=\inf\{\ell_G(xh)\ | \ h\in H\}=\ell_{G/H}(xH).
\end{align*}

Finally, we have $0\leq\ell_{G/H}(H)=\inf\{\ell_G(h)\ | \ h\in H\}\leq \ell_G(e_G)=0.$ 
\end{proof}

\begin{lemma}
Let $\phi: G\to H$ be a homomorphism between two groups and $\ell_G$ a length function on $G$. We can define the push-forward $\phi_\bullet(\ell_G): \phi(G)\to \RR$ as follows:
\begin{equation*}
(\phi_\bullet(\ell_G))(h):=\inf\{\ell_G(g)\ | \ g \in G, \ \phi(g)=h\}, \quad \forall h\in \phi(G).
\end{equation*}
The push-forward is a length function on $\phi(G).$
\end{lemma}
\begin{proof}
Under the natural isomorphism $\phi(G)\simeq G/\ke \phi,$ the function $\phi_\bullet(\ell_G)$ coincides with the function $\ell_{G/\ke \phi}$ defined above.
\end{proof}

\begin{remark}
A co-weighted homomorphism $\phi: G\to H$ between two charged groups is a co-isometry if and only if $|\omega_H|=\phi_\bullet(|\omega_G|).$
\end{remark}

\begin{lemma}
There is a category $\Gf_c$ whose objects are groups and whose morphisms are epimorphisms.

There is a category whose objects are charged groups and whose morphisms are co-isometric homomorphisms.
\end{lemma}
\begin{proof}
The composition of epimorphisms is another epimorphism and the composition of co-isometric homomorphisms 
is a co-isometric homomorphism. The identity map of every group is a co-isometric homomorphism (hence epimorphism) 
that plays the role of the identity in the category.
\end{proof}

\begin{corollary}
There is a category $\Hf_c$ whose objects are Hilbert spaces and whose morphisms are co-isometries.
\end{corollary}

\begin{definition}
A \emph{covariant relator} from the category $\Af$ to the category $\Bf$ is a pair $(\Rs_{\Ob},\Rs_{\Mor})$ of relations, $\Rs_{\Ob}\subset \Ob_\Af\times\Ob_\Bf$ between objects and $\Rs_{\Mor}\subset\Mor_\Af\times\Mor_\Bf$ between morphisms, such that: 
\begin{equation*}
(A,B)\in \Rs_{\Ob} \imp (\iota_A,\iota_B)\in \Rs_{\Mor}
\end{equation*}
and, whenever $\alpha_1, \alpha_2$ are composable morphisms in $\Af$ and whenever $\beta_1, \beta_2$ are composable morphisms in $\Bf:$  
\begin{equation*}
(\alpha_2,\beta_2), (\alpha_1,\beta_1) \in \Rs_{\Mor} 
\imp (\alpha_2\circ \alpha_1,\beta_2\circ \beta_1)\in \Rs_{\Mor}. 
\end{equation*}
\end{definition}
\begin{remark}
A covariant functor is a covariant relator such that both $\Rs_{\Ob}$ and $\Rs_{\Mor}$
are functions. Contravariant relators are defined in a similar way interchanging the order of compositions. 
\end{remark}

\begin{proposition}\label{pr: g-h2}
There is a contravariant relator $\H$ from the category $\Gf_c$ of groups with epimorphisms to the category $\Hf_c$ of Hilbert spaces with isometries.
\end{proposition}
\begin{proof}
$\H_{\Ob}$ is the function that to every object $G$ in $\Gf_c$ associates the Hilbert space $\H_G\in \Hf_c.$

We now define the relation $\H_{\Mor}.$
As we already know from proposition~\ref{pr: g-a}, every homomorphism $\phi: G\to H$ is associated to a linear map $\A_\phi: \A_G\to\A_H$ of pre-Hilbert spaces.

If $\phi$ is an epimorphism, the linear map $\A_\phi$ is continuous if and only if 
$\ke \phi$ is a finite subgroup of $G:$   
\begin{equation}\label{eq: in}
\begin{aligned}
\|\A_\phi(f)\|_H^2 &=\Big\|\A_\phi\Big(\sum_{x\in G}f(x)\delta_x^G\Big)\Big\|_H^2=\Big\|\sum_{x\in G}f(x)\delta_{\phi(x)}^H\Big\|_H^2\\
&= \Big\|\sum_{y\in H}\Big( \sum_{x\in \phi^{-1}(y)} f(x)\Big)\delta_y^H \Big\|^2_H= 
\sum_{y\in H}\Big|\sum_{x\in \phi^{-1}(y)}  f(x)\Big|^2  \\ 
& \leq \sum_{y\in H}\Big(\sum_{x\in \phi^{-1}(y)}|f(x)|\Big)^2 \leq 
\sum_{y\in H}\sum_{x\in \phi^{-1}(y)}\car(\ke \phi)|f(x)|^2  \\
&= \sum_{x\in G}\car(\ke \phi)|f(x)|^2 =\car(\ke \phi)\|f\|^2_G.
\end{aligned}
\end{equation}
It follows that in general the operator $\A_\phi$ is an unbounded operator from the Hilbert space $\H_G$ to the Hilbert space $\H_H.$ 

The operator $\A_\phi$ is densely defined because its domain contains the dense subspace $\A_G\subset \H_G.$ Hence there exists an adjoint operator $\A_\phi^*$ that, in the case of finite $\ke \phi,$ coincides with the ``pull-back'' operator $f\mapsto f\circ \phi,$ for  $f\in \A_H.$ 
Unfortunately, when $\ke \phi$ is not finite, $\A_\phi$ is not a closable operator.

Let us now denote by $\P$ the set of linear isometric operators $\K\subset \H_G\times\H_H$ such that $\K\subset \A_\phi.$
The family $\P$ is an inductive partially ordered set and as such, by Zorn's lemma, it admits a maximal element. Every maximal operator $\K$ in the family $\P$ is necessarily surjective and so its adjoint $\K^*:\H_H\to\H_G$ is an isometry and $\K^{**}=\cj{\K}$ is a partial isometry with range $\H_H.$ 

Every maximal partial isometry $\K\in \P$ has a closure that is the adjoint operator
$\H_\psi^*$ of an isometric operator $\H_\psi,$ where $\psi: H\to G$ is a monomorphism that is right inverse to the epimorphism $\phi: G\to H.$ 

We define a contravariant relator $\H_{\Mor}$ on morphisms by saying that $(\phi,\K^*)\in \H_{\Mor}$ if and only if $\K$ is a maximal isometry in $\A_\phi.$ 

$\H$ will be a contravariant functor if and only if the set of maximal partial isometries in $\A_\phi$ has cardinality one, which is equivalent to the fact that there exists only one ``splitting homomorphism'' for $\phi$ i.e.~there exists a unique $\psi:H\to G$ such that $\phi\circ\psi=\iota_H.$

The same considerations can be applied to the full subcategory of (proper) weighted groups. 
\end{proof}

\begin{remark}
The relation between the kernel of $\phi$ and the kernel of $\A_\phi$ is given by:
\begin{equation*}
\ke(\phi)=G\cap[\delta_e^G+\ke(\A_\phi)], \quad \ke(\A_\phi)=\delta_e^G + \spa(\ke(\phi)).
\end{equation*}
\end{remark}

\subsection{The Functor: Coisometric Case.}

\begin{theorem}
There exists a contravariant relator $\Rg$ from the category $\Gf^\omega_{c}$ of proper Dirac-weighted groups with 
co-weighted homomorphisms to the category of spectral triples $\Sf$, 
that to every $(G,\omega_G)$ associates $(\A_G,\H_G,D_{\omega_G})$.
\end{theorem}
\begin{proof}
The relator on objects $\Rg_{\Ob}$ coincides with the functor defined in 
theorem~\ref{th: im}. 

We want to see that, on morphisms, the relator $\Rg_{\Mor}$ is defined in the same way as in proposition~\ref{pr: g-h2} i.e.~$\Rg_{\Mor}$ associates to every splitting weighted epimorphism $\phi: G\to H$ the family of pairs $(\A_\psi, \H_\psi),$ where $\psi$ is any weighted (mono)morphism $\psi:H\to G$ such that $\phi\circ \psi=\iota_H$.

Let $\phi: G\to H$ be a co-weighted epimorphism of proper weighted groups. 
For sure (see proposition~\ref{pr: g-a}) we have that $\A_\psi: \A_G\to \A_H$ is an involutive unital homomorphism of the group algebras. 

If the homomorphism $\phi:G\to H$ admits ``right inverses'' i.e.~if there exist weighted morphisms $\psi: H\to G$ such that $\phi\circ\psi=\iota_H$, from proposition~\ref{pr: g-h2} we know that, for any such ``right inverse'' $\psi$, the function $\H_\psi:\H_G\to \H_H$ is an isometry of Hilbert spaces. 
From the same proposition~\ref{pr: g-h2} we also know that the pair $(\Rg_{\Ob},\Rg_{\Mor})$ where the second relation is given by $\Rg_{\Mor}:=\{(\A_\psi,\H_\psi)\ | \ \phi\circ\psi=\iota_H, \ \phi\in \Mor_{\Gf^\omega_c} \}$ is a contravariant relator.

From theorem~\ref{th: im}, $(\A_\psi,\H_\psi)$ is a morphism in the category of spectral triples $\Sf$, i.e.~for all $f\in \A_H$ and for all $\xi\in \A_H:$
\begin{gather*}
\H_\psi \circ \pi_H(f) (\xi)= \pi_G(\A_\psi(f))\circ \H_\psi (\xi)\\
D_{\omega_G}\circ \H_\psi(\xi)=\H_\psi\circ D_{\omega_H}(\xi).
\end{gather*}
\end{proof}

\begin{remark}
The relator becomes a functor in case that it is possible to select canonically a splitting of the co-isometry (for example in the case of Hilbert spaces). 
\end{remark}

\section{Conclusion and Further Remarks.}

In this work we have proposed a definition of morphism for spectral triples (and their real and even variants). 

Some remarks on further generalizations are in order.
We have presented here the most elementary instructive example of functorial relations between our proposed category $\Sf$ of spectral triples and other categories: in this specific case the categories $\Gf^\omega_i$ of proper Dirac-weighted groups with monomorphisms and $\Gf^\omega_c$ of proper Dirac-weighted groups with co-weighted homomorphisms.
Other examples involving categories of Riemannian manifolds equipped with a spin$^c$ structure will be dealt with in~\cite{BCL2}.

Other alternative variants of our definition of morphism of spectral triples are worth investigating. 
For example we might substitute the ``strong'' requirement of commutation of the Dirac operators with the Hilbert space maps with some milder property like\footnote{Similar properties in the case of unitary maps between Hilbert spaces have been recently suggested by M.~Paschke and R.~Verch~\cite[page~10]{PV}.}: $\H_\phi\circ[D_G, \pi_G(f) ]=[D_H,\pi_H(\A_\phi(f))]\circ \H_\phi.$

\medskip

As regards the specific examples of functorial relations described here, several immediate generalizations and comments come to mind. 
Among them, we mention:
\begin{itemize}
\item[]
Most of the facts presented here for the category of weighted groups can be rephrased for  the category of ``weighted'' small categories considering the ``convolution algebra'' of a  small category in place of the group algebra. 
\item[]
The notions of weight and charge on a group can be further generalized by considering functions $\omega:G\times G\to \CC$ having properties formally similar to those of Hermitian forms and inner products. The Dirac operators $D_{\omega}$ associated to these functions $\omega$ include, in the case of finite groups, all available Dirac operators according to the classification of finite spectral triples (see, for example, T.~Krajewski~\cite{K}).
\item[]
The only possible choice of Dirac operator $D_{\omega_G}$ on a weighted group $(G,\omega_G)$ that is fully compatible with the requirements of a real zero dimensional spectral triple, where the real structure $J$ is the one obtained by Tomita-Takesaki Modular theory (from the cyclic separating vector $\delta^G_{e_G}\in \H_G$), is $D_{\omega_G}=0$.
This fact asks for some investigation on the mutual relationship between modular theory and non-commutative geometry. 
We hope to discuss this point elsewhere~\cite{BCL1}.
\end{itemize}

\bigskip

{\bf Acknowledgments.}
We acknowledge the support provided by the Thai Research Fund through the ``Career Development Grant'' n.~RSA4580030: ``Modular Spectral Triples in Non-commutative Geometry and Physics''.

We thank M.~Marcolli for some enlightening discussion and E.~B\'edos for some useful comments on a preliminary draft of this work.

We also thank A.~Rennie for providing a draft of his papers~\cite{Re1,Re2} prior to publication.


\begin{thebibliography}{XXXXX}

\bibitem[Ba]{Ba}
J.~Baez, Categories, Quantization, and Much More, \\
\hlink{http://math.ucr.edu/home/baez/categories.html}, 26 September 2004.

\bibitem[BD]{BD}
J.~Baez, J.~Dolan, Categorification, Contemp.~Math. \emph{230}, 1-36, (1998), 
\hlink{http://xxx.lanl.gov/math/9802029}, 05 February 1998. 

\bibitem[BCL1]{BCL1}
P.~Bertozzini, R.~Conti, W.~Lewkeeratiyutkul, Modular Spectral Triples, in progress.

\bibitem[BCL2]{BCL2}
P.~Bertozzini, R.~Conti, W.~Lewkeeratiyutkul, Non-commutative Totally Geodesic Submanifolds and Quotient Manifolds, in preparation.

\bibitem[BCL3]{BCL3}
P.~Bertozzini, R.~Conti, W.~Lewkeeratiyutkul, Categories of Spectral Triples and Morita Equivalence, in progress.

\bibitem[C1]{C1}
A.~Connes, Compact Metric Spaces, Fredholm Modules and Hyperfiniteness, Ergodic Theory Dynam.~Systems \emph{9(2)}, 207-220 (1989).

\bibitem[C2]{C2}
A.~Connes, Noncommutative Geometry, Academic Press (1994).

\bibitem[C3]{C3}
A.~Connes, Noncommutative Geometry and Reality, J.~Math.~Phys. \emph{36}, No.~11, 6194-6231 (1995). 

\bibitem[C4]{C4}
A.~Connes, Gravity coupled with Matter and the Foundations of Noncommutative Geometry, Commun.~Math.~Phys. \emph{182}, 155-176 (1996).


\bibitem[C5]{C5}
A.~Connes, Noncommutative Geometry Year 2000, GAFA special volume 2000 (2001), \hlink{http://xxx.lanl.gov/math.QA/0011193}.

\bibitem[FGV]{FGV} 
J.~M.~Gracia-Bondia, H.~Figueroa, J.~C.~Varilly, Elements of Noncommutative Geometry, Birkh\"auser (2000).

\bibitem[GGISV]{GGISV}
V.~Gayral, J.~M.~Gracia-Bondia, B.~Iochum, T.~ Sch\"uker, J.~C.~Varilly, Moyal Planes are Spectral Triples, Commun.~Math.~Phys. \emph{246}, no.~3, 569-623 (2004),
\\ \hlink{http://xxx.lanl.gov/hep-th/0307241}, 24 July 2003.

\bibitem[GLMV]{GLMV}
J.~M.~Gracia-Bondia, F.~Lizzi, G.~Marmo, P.~Vitale, Infinitely Many Star Products to Play with, J.~High~Energy~Phys. \emph{4}, n.~26 (2002), 
\hlink{http://xxx.lanl.gov/hep-th/0112092}, 19 December 2001.

\bibitem[I]{I}
L.~Ionescu, On Categorification, \\
\hlink{http://xxx.lanl.gov/math.CT/9906038}, 06 June 1999.
% Categorification and group extensions. Appl. Categ. Structures  10  (2002), no. 1, 35--47

\bibitem[KR1]{KR1}
M.~Kontsevich, A.~Rosenberg, Noncommutative Smooth Spaces, in:
The Gel'fand Mathematical Seminars 1996-1999, 85-108, Birkh\"auser (2000), 
\hlink{http://xxx.lanl.gov/math.AG/9812158}, \\
30 December 1998.

\bibitem[KR2]{KR2}
M.~Kontsevich, A.~Rosenberg, Noncommutative Spaces,\\ 
preprint~MPIM2004-35, Max Planck Institut f\"ur Mathematik (2004).

\bibitem[K]{K}
T.~Krajewski, Classification of Finite Spectral Triples, \\ 
J.~Geom.~Phys. \emph{28}, 1-30 (1998), \\  
\hlink{http://xxx.lanl.gov/hep-th/9701081}, 20 January 1997.

\bibitem[Ma]{Ma}
Y.~Manin, Real Multiplication and Noncommutative Geometry, in:
The Legacy of Niels Henrik Abel, 685-727, Springer (2004). \\ 
\hlink{http://xxx.lanl.gov/math.AG/0202109}, 12 February 2002.

\bibitem[PV]{PV}
M.~Paschke, R.~Verch, Local Covariant Quantum Field Theory over Spectral Geometries, 
Class.~Quantum~Grav. \emph{21}, n.~23, 5299-5316 (2004), 
\hlink{http://xxx.lanl.gov/gr-qc/0405057}, 11 May 2004.

\bibitem[Re1]{Re1}
A.~Rennie, Smoothness and Locality for Nonunital Spectral Triples, $K$-Theory \emph{28}, n.~2, 127-165 (2003). 
% preprint, 11 August 2003. 

\bibitem[Re2]{Re2}
A.~Rennie, Summability for Nonunital Spectral Triples, $K$-Theory \emph{31}, n.~1, 71-100  (2004). 
%preprint, 11 August 2003. 

\bibitem[Ri1]{Ri1}
M.~A.~Rieffel, Metrics on States from Action of Compact Groups, Doc.~Math. \emph{3}, 215-229 (1998), \\
\hlink{http://xxx.lanl.gov/math/9807084}, 04 January 1999.

\bibitem[Ri2]{Ri2}
M.~A.~Rieffel, Group C$^*$-algebras as Compact Quantum Metric Spaces, 
Doc.~Math. \emph{7}, 605-651 (2002), \\
\hlink{http://xxx.lanl.gov/math.OA/0205195}, version~3, 21 November 2002.

\bibitem[R]{R}
A.~Rosenberg, Noncommutative Spaces and Schemes, \\
preprint~MPIM1999-84, 
Max Planck Institut f\"ur Mathematik (1999).

\end{thebibliography}
\end{document}